\documentclass[a4paper,12pt,reqno]{amsart} 

\usepackage[margin=2cm]{geometry}
\usepackage[utf8]{inputenc}
\usepackage[T1]{fontenc}
\usepackage{lmodern}
\usepackage[french,english]{babel}
\usepackage{amsthm, amssymb, amsmath, amsfonts, mathrsfs, dsfont, esint, textcomp, bm}
\usepackage{mathtools}
\usepackage{tikz}
\usepackage{enumitem, imakeidx}
\usetikzlibrary{patterns}

\usepackage{amsthm, amssymb, amsmath, amsfonts, mathrsfs, dsfont, esint}
\usepackage[colorlinks=true, pdfstartview=FitV, linkcolor=blue, citecolor=blue, urlcolor=blue,pagebackref=false]{hyperref}
\theoremstyle{remark}

\theoremstyle{definition}

\renewcommand{\leq}{\leqslant}
\renewcommand{\geq}{\geqslant}

\renewcommand{\subset}{\subseteq}

\newcommand{\F}{\mathcal{F}}

\newcommand{\N}{\mathbb{N}}

\newcommand{\R}{\mathbb{R}}

\newcommand{\eps}{\varepsilon}
\renewcommand{\d}{{\mathrm{d}}}

\newcommand{\Rmm}{\R^{m\times m}}
\newcommand{\Rm}{\R^{m}}
\newcommand{\pcap}{\text{{\small p-}cap}}
\newcommand{\Rd}{\mathbb{R}^d}

\newcommand{\ignore}[1]{}
\newcommand{\oh}{{\textstyle\frac{1}{2}}}
\newcommand{\M}{\mathcal{M}}

\newtheorem{theorem}{Theorem}

\newtheorem{corollary}{Corollary}
\providecommand{\dx}{\, \mathrm{d} x}
\providecommand{\dy}{\, \mathrm{d} y}

\title{On the existence of the Green function for elliptic systems in divergence form}
\author{Arianna Giunti, Felix Otto}

\begin{document}
\maketitle

\begin{abstract}
We study the existence of the Green function for an elliptic system in divergence form $-\nabla\cdot a\nabla$ in $\Rd$, with $d>2$. The tensor field 
$a=a(x)$ is only assumed to be bounded and $\lambda$-coercive.  For almost every point $y \in \Rd$, the existence of a Green's function $G(\cdot, y)$ centered in $y$ has been proven in \cite{ConlonGiuntiOtto}. In this paper we show that the set of points $y \in \Rd$ for which $G( \cdot, y)$ does not exist has zero $p$-capacity, for an exponent $p >2$ 
depending only on the dimension $d$ and the ellipticity ratio of $a$. 
\bigskip


\medskip


\end{abstract}

This paper is an extension of \cite{ConlonGiuntiOtto} and further investigates the existence of a Green's function for the second-order 
elliptic operator $-\nabla \cdot a \nabla$ in  $\Rd$, with $d> 2$. We focus on the case of systems of $m$ equations, namely when $a$ is a measurable tensor field 
$a : \Rd \rightarrow \mathcal{L}(\mathbb{R}^{m \times d} ; \mathbb{R}^{m \times d})$, with $m$ being any positive integer. 
We stress that in this paper we do only assume that $a$ is bounded and $\lambda$-coercive, i.e.\ that there exists $\lambda >0$ such that
\begin{equation}\label{a.assumptions}
\begin{aligned}
& \forall \zeta \in C^\infty_0( \Rd ; \R^m)\ \ \ \int \nabla \zeta(x) \cdot a(x) \nabla \zeta(x) \dx \geq \lambda \int |\nabla \zeta(x)|^2 \dx,\\ 
& \forall x\in \mathbb{R}^d,\ \forall \ \xi\in \R^{m\times d} \ \ \ |a(x)\xi|\leq |\xi| .
\end{aligned}
\end{equation}

\medskip

In \cite{ConlonGiuntiOtto}, J. Conlon and the authors show that a Green's function $G(\cdot , y)$ centered in $y$ exists for every coefficient field $a$ satisfying \eqref{a.assumptions} and for (Lebesgue-)almost every point  $y\in \Rd$. In this paper, we improve this result by showing that the \textit{exceptional set} $\Sigma$ of points $y\in \Rd$ 
for which $G(a; \cdot, y)$ does not exist has p-capacity zero, for an exponent $p > 2$ depending only on the dimension $d$ and the ellipticity ratio $\lambda$. This, in particular, implies that for every coefficient field $a$ that is $\lambda$-coercive and bounded, the Hausdorff dimension of $\Sigma$ is strictly smaller than $d-2$ \cite{EG}[Theorem 4.17]. 

\medskip

The result of \cite{ConlonGiuntiOtto} crucially relies on the idea of studying the Green function as a map $G = G(\cdot, \cdot)$ in both variables $x, y \in \Rd$. This yields optimal estimates for the $L^2$-norm in $y$ and $x$ of $G$, $\nabla_x G$ and $\nabla_x \nabla_y G$ both away from the diagonal $\{ x = y \}$ and close to it. By the standard properties of Lebesgue-integrable functions, these estimates allow to give a pointwise meaning in $y$ to $G( \cdot , y)$, up to a set of Lebesgue-measure zero. The main idea behind the result of this paper is to exploit the integrability of the mixed derivatives $\nabla_y \nabla_x G$ and extend the set of Lebesgue points $y$ where $G( \cdot , y)$ is well-defined up to the set $\Sigma$ having zero $p$-capacity.

\medskip

We remark that in the case of elliptic systems the set $\Sigma$ may indeed be non-trivial. There are, indeed, coefficient fields $a$ satisfying \eqref{a.assumptions} for which one may construct \textit{unbounded} $a$-harmonic vector fields. From this, and by means of representation formulas, it follows that the points where such vector fields are unbounded \textit{cannot} be Lebesgue points for $G( \cdot, y)$. A classical example of a discontinuous $a$-harmonic vector field is due to E. De Giorgi \cite{DeGiorgiCounterexample}: For any dimension $d> 2$, the vector field $u : \Rd \to \Rd$
\begin{align}\label{De.Giorgi.u}
u(x)=\frac{x}{|x|^{\gamma}}, \text{\ \ \ $\gamma:=\frac{d}{2}\biggl(1-\frac{1}{\sqrt{(2d-2)^2+1}}\biggr)> 1$}
\end{align}
solves $-\nabla \cdot a_0 \nabla u = 0$  in $\Rd$, with $a_0$ satisfying \eqref{a.assumptions} and being smooth everywhere outside of the origin. We remark that the coefficient $a_0$ is not only $\lambda$-coercive as in \eqref{a.assumptions}, but also  {\it strongly elliptic}: For almost every $x \in \R^d$ and every matrix $\xi \in \R^{d \times d}$, it satisfies $\xi \cdot a_0(x) \xi \geq \lambda |\xi|^2$, with $\lambda$ depending on $d$.

\medskip

In the case $d=3$, the previous example implies that the exceptional set $\Sigma$ for $a_0$ contains at least the origin. For higher dimensions $d \geq 3$, the trivial extension of the vector field $u$ for $d=3$ is itself $\bar a_0$-harmonic if 
$$
\bar a_0:=\begin{pmatrix}
a_0 & 0 \\
0 & \mathop{I}
\end{pmatrix}.$$
This implies, in particular, that $\Sigma$ for $\bar a_0$ has Hausdorff dimension at least $d-3$. 
 
\smallskip

The previous counterexample also implies that for (locally) $a$-harmonic vector fields one may only aim at statements on their \textit{partial regularity} as, for instance, their continuity outside of a singular set. We remark that there exist examples of discontinuous $a$-harmonic vector fields with discontinuity much larger than \eqref{De.Giorgi.u}: We refer, for instance, to the paper by J. Soucek \cite{Soucek}, which exhibits an $a$-harmonic vector field discontinuous on a dense countable set, and the one by O. John, J. Mal\'y and J. Star\'a \cite{JohnMalyStara}, in which, for every countable union of closed sets (i.e. $F_\sigma$-set), an $a$-harmonic vector field discontinuous there is constructed. 

\smallskip

Without using the equation, the fact that $a$-harmonic functions are locally in $H^1$ immediately implies that they are \textit{2-quasicontinuous}. This means that there exist sets, of arbitrarily small $2$-capacity, outside of which the function considered is continuous \cite{EG}[Definition 4.11]. This argument is oblivious to the difference between scalar and vectorial functions. Using the equation and appealing to Meyers's \cite{MeyersLp} or Gehring's \cite{giaquinta} estimates, this notion of continuity may be upgraded from $2$-quasicontinuity to p-quasicontinuity, for an exponent $p>2$. The result of this paper provides an analogous statement for the solution operator for $-\nabla \cdot a \nabla$. By means of representation formulas, indeed, we prove that for any family $\mathcal{F}$ of locally $a$-harmonic functions that are uniformly bounded in the $H^1_{loc}$-norm, there exist \textit{common} sets of arbitrarily small $p$-capacity outside of which $\mathcal{F}$ is equicontinuous (see Corollary \ref{cor1}).  These sets are universal in the sense that they depend only on the coefficient $a$ and on the dimension $d$, but not on the family  $\mathcal{F}$. 

 \bigskip

{\bf Notation and previous results.} For notational convenience, as in \cite{ConlonGiuntiOtto} we assume
that $a$ is symmetric, i.e. that for almost every $x\in \Rd$, the tensor $a(x) \in \mathcal{L}(\mathbb{R}^{m \times d} ; \mathbb{R}^{m \times d})$ is symmetric. Throughout this paper the expression ``almost every'' is meant with respect to the Lebesgue measure and all the PDEs considered are assumed to hold only in the distributional sense. 

\smallskip
 
We denote by $W^{1,p}(\Rd, \R^m)$, $p \geq 1$ the Sobolev spaces of functions in
$\Rd$ taking values in $\R^m$; if $m =1$, we use the usual notation $W^{1,p}(\Rd)$. The same criteria are employed for all the
other standard functions spaces used in the paper. If we represent the elements of the product space $\R^d \times \R^d$ by $(x,y)$, we use the notation 
$W^{1,q}_x(\Rd)$ to specify in the lower index the differentiation and integration variable. Similarly, we write $\nabla_x, \nabla_{y}$ or $\nabla_{x,y}$ when the gradient is taken with respect to $x$, $y$ or both variables $(x,y)$, respectively. We write $\gtrsim$ and $\lesssim$ for $\geq C$  and $\leq C$ with a constant $C$ depending only on the dimension $d$, the 
ellipticity ratio $\lambda$ and the dimension of the target space $m$.

\medskip

{For an open (not necessarily bounded) set $D \subset \Rd$, we may define the space
$$Y^{1,2}(D ; \Rm) := \bigl\{ u \in L^{\frac{2d}{d-2}}(D ; \Rm) \, : \, \nabla u \in L^2(D ; \R^{m \times d}) \bigr\},$$
and equip it with the norm $\| u\|_{Y^{1,2}(D ; \Rm)}:= \| u \|_{L^{\frac{2d}{d-2}}(D; \Rm)} + \| \nabla u \|_{L^2(D ; \R^{m \times d})}$.} 
The main theorem of \cite[Theorem 1]{ConlonGiuntiOtto} provides the existence of a map 
\begin{align}
\begin{aligned}\label{G.map}
G& : \Rd \times \Rd \rightarrow \R^{m \times m}
\end{aligned}
\end{align} 
such that for all $ 1\leq q < \frac{d}{d-1}$
\begin{align}\label{G.L^q}
G \in W^{1, q}_{loc}(\Rd \times \Rd ; \R^{m \times m})
\end{align}
and for almost every $y \in \Rd$, the tensor field $G( \cdot , y ) \in Y_x^{1,2}(\{ |x-y| > 1 \} ; \R^{m\times m})$ and satisfies
\begin{equation}\label{equation.Green}
-\nabla_x \cdot a(\cdot) \nabla_x G(\cdot , y) = \delta( \cdot - y) \ \ \ \ \ \ \text{in $\Rd$.}
\end{equation}
Furthermore, for every $R > 0$, $z\in \mathbb{R}^d$ and $\alpha > \frac d 2 -1$ it holds 
\begin{align}
\int_{|y-z|< R}\int_{|x-z|< R} |x-y|^{2\alpha} |\nabla_{x,y} G(x, y)|^2  dx \dy  &\lesssim_\alpha R^{2+ 2\alpha},  \label{weighted.map} \\
\int_{|y-z|< R}\int_{|x-z|> 2R} |\nabla_x G(x, y)|^2 dx \dy &\lesssim R^2, \label{grad.decay}\\
\int_{|y-z|< R}\int_{|x-z|> 2R} |\nabla_y\nabla_x G(x, y)|^2 dx \dy &\lesssim 1 \label{mixed.decay}.
\end{align}
Here, the notation $\lesssim_\alpha$ means that the constant depends also on the exponent $\alpha$.

\section{Main result}
\begin{theorem}\label{t.main}
Recall $d > 2$. Let $a$ be symmetric and satisfy assumptions \eqref{a.assumptions}. Let \nolinebreak{$G : \mathbb{R}^d \times \mathbb{R}^d \rightarrow \mathbb{R}^{m \times m}$} be the Green function in the sense of \cite{ConlonGiuntiOtto}
and constructed there. Then, there exists a (measurable) set $\Sigma= \Sigma(a) \subset \Rd$ with the following properties:
\begin{itemize}
\item[(a)] There exists an exponent $p= p(d, \lambda) > 2$ such that
$$\pcap \bigl( \Sigma \bigr) = 0. $$

\item[(b)]  There exists an exponent $q=q(d, \lambda) >1$ such that for every $y \in \Rd \setminus \Sigma$ and for all $R >0$
\[  \fint_{|z- y| < \delta }\hspace{-0.4cm} G( \cdot ,  z) \, d z \rightharpoonup G^*( \cdot , y) \ \ \ \ \text{as $\delta \downarrow 0$} \]
in $W^{1, q}_{ x}(\{ |x-y| < R \} ; \Rmm) \cap Y^{1,2}_x(\{ | x- y | > R \}; \Rmm)$. \vspace{0.1cm}

\item[(c)] The representative $G^*( \cdot , y)$ solves equation \eqref{equation.Green} for every $y \in \Rd \setminus \Sigma$. Furthermore, there exists an exponent $\alpha < \frac d 2$ such that
\[ \int \min\{ |x-y|^{2\alpha} , 1 \} |\nabla_x G^*( x, y)|^2 \dx < + \infty.\]
\end{itemize}
\end{theorem}
In addition, as a corollary we have:
\begin{corollary}\label{cor1}
Let  the coefficient field $a$ and the exponent $p > 2$ be as in Theorem \ref{t.main}. Consider
\begin{align*}
\mathcal{F}:=\biggl\{ u \, \colon \, \ \text{u is $a$-harmonic in $\{ |x| < 4\}$, $\int_{|x| < 4} |\nabla u|^2 \leq 1$} \ \ \biggr\}.
\end{align*}
Then for every $\eps > 0$ there exists an open set $U^\eps \subset \{|x| < 1\}$ having 
$$ \pcap(U^{\varepsilon}) < \varepsilon $$
such that $\F$ is equicontinuous in $\{|x| < 1\}\backslash U^\eps$ and the modulus of continuity is uniform in $\{|x| < 1\}\backslash U^\eps$.
\end{corollary}

\section{Proofs}
For the sake of simplicity, throughout this section we use a scalar notation and language by pretending that $a$ is a field of $d\times d$ matrices and that the Green function is scalar. For a detailed
discussion about this abuse of notation, we refer to \cite[Section 2]{ConlonGiuntiOtto}. Moreover, for a function $f \in L^1_{loc}(\Rd)$, we introduce the notation
 \begin{align}\label{max.function}
  \M f(y) &:= \sup_{\delta > 0} \fint_{|\tilde y-y|< \delta} | f(\tilde y)| \, d\tilde y,\\
  \M_0 f(y) &:= \limsup_{\delta \downarrow 0}\fint_{|\tilde y -y|< \delta} |f(\tilde y)| \, d\tilde y.\label{limsup}
 \end{align}

\medskip

Before giving the proof of Theorem \ref{t.main}, we recall some further properties of $G(\cdot, \cdot)$ obtained in \cite[Section 2]{ConlonGiuntiOtto} which will be used in our proofs:
\begin{itemize}

\smallskip

\item For every $R >0$, $ z\in \mathbb{R}^{d}$ and almost every $y \in \Rd$ such that $|y- z| < R$,
\begin{align}\label{mixed.equation}
 -\nabla_x \cdot a(\cdot) \nabla_x \nabla_{y} G(\cdot , y) = 0 \mbox{\ \ \hspace{0.8cm}   in $\{|x- z| > 2R \}$.}
\end{align}
and $\nabla_y G(\cdot, y) \in Y^{1, 2}_x(\{|x- z| > 2R \}, \Rd)$.

\smallskip

\item For almost every $(x, y) \in \Rd \times \R^d$
\begin{align}\label{symmetry.G}
 G( x, y) = G(y, x).
\end{align}

\smallskip

\item  {For every $g \in L^2(\Rd; \Rd)$ and every $f \in L^2(\Rd)$ having compact support, the  solution $u \in Y^{1,2}(\Rd)$ to 
$$
-\nabla \cdot a \nabla u = \nabla \cdot g +f \ \ \text{in $\Rd$}
$$ 
may be written as the identity (up to a set of Lebesgue-measure zero)
\begin{align}\label{representation.formula}
u(\cdot) = - \int  g(y) \cdot \nabla_y G( \cdot, y) \dy + \int  f(y) G(\cdot, y) \dy.
\end{align}
}
\end{itemize}

\noindent {\it Proof of Theorem \ref{t.main}.\,} We divide the proof into five steps: In Step 1 we give a formulation of the standard Gehring's estimate tailored to our needs. It allows to upgrade estimate \eqref{mixed.decay} into an $L^2$-estimate in $x$ and $L^{\bar p}$ in $y$, for the Gehring exponent $\bar p > 2$. Steps 2-4 contain the main capacitary estimates for the exceptional set $\Sigma$, which is closely related to the set of points $y \in \Rd$ where $G( \cdot , y)$ and $\nabla_x G( \cdot, y)$ have infinite $W^{1,q}_{x}(\{|x-y| < 1 \})$- or $Y^{1,2}_x( \{ |x- y | > 1 \})$-norms. These estimates on the capacity of $\Sigma$ crucially rely on the upgraded version of \eqref{mixed.decay} and are combined with a maximal function estimate for Sobolev functions. Finally, in Step 5 we argue how to construct the representative $G^*( \cdot, y)$ away from the singularity set $\Sigma$.   

\bigskip

{\bf Step 1. Gehring's estimate.} \, Let  $u \in H^1(\{ |x| < 2R\})$  be a solution to
\begin{align}\label{equation.Gehring}
- \nabla \cdot a \nabla u = \nabla \cdot g \ \ \ \ \ \ \ \text{ in $\{ |x| < 2R \}$.}
\end{align}
Then, there exists an exponent $\bar p=\bar p(d, \lambda) > 2$ such that
\begin{align}\label{Gehring}
\bigl( \fint_{|x| < R} |\nabla u(x)|^{\bar p} \dx \bigr)^{\frac{1 }{\bar p}} \lesssim \bigl( \fint_{|x| < 2R} |\nabla u(x)|^2 \dx \bigr)^{\frac 1 2} + \bigl( \fint_{|x| < 2R} | g(x)|^{\bar p} \dx \bigr)^{\frac {1}{ \bar p}}.
\end{align}
This is a standard result in elliptic regularity theory and we refer to \cite[Chapter V, Theorem 2.1]{giaquinta} for its proof.\footnote{ In \cite{giaquinta}, the coefficients are assumed to be strongly elliptic (cf. \cite[Chapter V, display below (0.1)]{giaquinta}). However, the argument only relies on 
 Caccioppoli's and Poincar\'e-Sobolev's inequality which hold true also if $a$ is only $\lambda$-coercive as in \eqref{a.assumptions}. Moreover, \cite{giaquinta}[Chapter V, Inequality (0.2)] corresponds to the standard case of $a$-harmonic functions; our case is an immediate adaptation of the Caccioppoli inequality in the case of solutions with a non-zero right-hand side as in \eqref{equation.Gehring}. } 
We pick a (smooth) cut-off function  $\eta$ for $\{|x| < R \}$ in $\{ |x| < 2R \}$. Since we may assume ${\bar p} \leq \frac{2d}{d-2}$,  the Poincar\'e-Sobolev inequality yields
\begin{align*}
\bigl( \fint_{|x| < 2R} |\eta(x) u(x)|^{\bar p} \dx \bigr)^{\frac{1}{\bar p}} \leq R \bigl( \fint_{|x| < 2R} |\nabla(\eta(x) u(x))|^2 \dx \bigr)^{\frac 1 2},
\end{align*}
it follows that
\begin{align}\label{Gehring2}
R^{-1}\bigl( \fint_{|x| < R}& | u(x)|^{\bar p }\dx \bigr)^{\frac {1}{\bar p}} + \bigl( \fint_{|x| < R} |\nabla u(x)|^{\bar p} \dx \bigr)^{\frac{ 1}{ \bar p}}\\
& \lesssim  R^{-1}\bigl( \fint_{|x| < 2R} | u(x)|^2 \dx \bigr)^{\frac 1 2}+  \bigl( \fint_{|x| < 2R} |\nabla u(x)|^2 \dx \bigr)^{\frac 1 2} + \bigl( \fint_{|x| < 2R} | g(x)|^{\bar p} \dx \bigr)^{\frac{1}{\bar p}}.\notag
\end{align}

\bigskip

 {\bf Step 2. Capacity estimates: First reduction.} \, Recall the definition \eqref{limsup} of $\mathcal{M}_0$ which we always think as acting on the $y$-variable. Let $p \geq 1$. We claim that if for an exponent $0 < \alpha=\alpha(d, \lambda) < \frac d 2 $ and every $R> 0$, $z \in \Rd$
 \begin{equation}
 \begin{aligned}
\pcap\biggl\{  \, |y -z|< 1 \, \colon \, \M_0 \biggl( \int_{|x-y| < R}\hspace{-0.2cm} |x-y|^{2\alpha} |G( x, y)|^2 \dx \biggr)^\oh = +\infty \biggr\} &= 0,\\\label{est.00b}
\pcap\biggl\{  \, |y -z|< 1 \, \colon \, \M_0 \biggl( \int_{|x-y| < R}\hspace{-0.2cm} |x-y|^{2\alpha} |\nabla_x G(x, y)|^2 \dx \biggr)^\oh = +\infty \biggr\} &= 0,
\end{aligned}
\end{equation}
then we may also find an exponent $q= q(d, \lambda) > 1$ such that for every $R> 0$
 \begin{align}
  \pcap\biggl\{ y \in \Rd \, \colon \, \M_{0} \|G(\cdot, y) \|_{W^{1,q}_x(\{ |x-y| < R \} ) } = +\infty \biggr\} = 0 \label{est.0}.
 \end{align}
Similarly, if for every $R> 0$ and $z \in \Rd$ it holds
\begin{equation}
\begin{aligned}\label{est.00}
\pcap\biggl\{  \, |y -z|< 1 \, \colon \, \M_0 \biggl( \int_{|x-y| > R} |G( x, y)|^{\frac{2d}{d-2}} \dx \biggr)^{\frac{d-2}{2d}} = +\infty \biggr\} &= 0,\\
 \pcap\biggl\{  \, |y -z|< 1 \, \colon \, \M_0 \biggl( \int_{|x-y| > R} |\nabla_x G( x, y)|^2 \dx \biggr)^\oh = +\infty \biggr\} &= 0,
\end{aligned}
\end{equation}
then for every $R> 0$
 \begin{align}
 \pcap\biggl\{ y \in \Rd \, \colon \, \M_{0} \|G(\cdot, y)\|_{Y^{1,2}_x(\{|x-y| > R \})} = +\infty \biggr\} = 0.\label{est.0a}
 \end{align}

\medskip

Since we may cover the whole space $\Rd$ with a countable number of unit balls, the subadditivity of the capacity and \eqref{est.00} immediately imply \eqref{est.0a}.
Analogously, the subadditivity of the capacity and \eqref{est.00b}  yield that for each $R > 0$
\begin{equation}\label{identities.d.sc}
 \begin{aligned}
  \pcap\biggl\{ y \in \Rd \, \colon \, \M_0 \bigl( \int_{|x-y| > R} |x-y|^{2\alpha} (|G(x, y)|^2 + |\nabla_x G(x, y)|^2 )\dx \bigr)^\oh = +\infty \biggr\} = 0.
 \end{aligned}
 \end{equation}
Since for $\alpha < \frac d 2$, H\"older's inequality implies that there exists $1 < q =q(\alpha) < 2$ such that for any $u$
\begin{align*}
 \int_{|x-y| < R} |u(x)|^q \dx &\lesssim_R \biggl( \int_{|x-y| < R} |x-y|^{2\alpha}|u(x)|^2  \dx \biggr)^{\frac q 2},
\end{align*}
estimate \eqref{est.0} is implied by this inequality together with identity \eqref{identities.d.sc} and the monotonicity of the capacity.

\bigskip

 {\bf Step 3. Capacity estimates: Second reduction.} \, Let $2 < p < \bar p$ be fixed, with $\bar p$ as in Step 1. We now argue that in order to prove $\eqref{est.00b} \& \eqref{est.00}$ with this choice of exponent $p$, it suffices to show that for every $R > 0$, $z \in \Rd$, and all $\lambda > 0$
\begin{equation}\label{est.C}
 \begin{aligned}
\pcap\biggl\{ |y-z| < \frac R 2 \, \colon \, \M_0\biggl( \int_{ |x-y| > 8R } |\nabla_{x}G( x, y)|^2 \dx \biggr)^{\frac {\bar p}{2p}}> \lambda \biggr\} &\lesssim \lambda^{-p} R^{-\frac{ d\bar p}{ 2} + d + (\bar p - p)}, \\
\pcap\biggl\{ |y-z| < \frac R 2 \, \colon \, \M_0\biggl( \int_{ |x-y| > 8R } |G( x, y)|^{\frac{2d}{d-2}} \dx \biggr)^{\frac {d-2}{2d}}> \lambda \biggr\} &\lesssim \lambda^{-p} R^{-\frac{ dp}{ 2} + d}.
\end{aligned}
\end{equation}

\medskip

Without loss of generality, we argue $\eqref{est.00b} \& \eqref{est.00}$ in the case $z=0$ and for $R=1$. We begin by observing that, since $p \leq \bar p$, definition \eqref{limsup} for $\mathcal{M}_0$ together with Jensen's inequality yield that for every $y \in \R^d$
$$
\M_0\biggl( \int_{ |x-y| > 8R } |\nabla_{x}G( x, y)|^2 \dx \biggr)^{\frac {\bar p}{2p}} \geq \biggl( \M_0\biggl( \int_{ |x-y| > 8R } |\nabla_{x}G( x, y)|^2 \dx \biggr)^{\frac {1}{2}}\biggr)^{\frac{\bar p}{p}}.
$$
This and the first inequality in \eqref{est.C} imply, after relabelling the parameter $\lambda^{\frac{p}{\bar p}}$ as $\lambda$, that for every $R > 0$ and $\lambda > 0$ it holds
\begin{align}\label{est.gradient.jensen}
\pcap\biggl\{ |y-z| < \frac R 2 \, \colon \, \M_0\biggl( \int_{ |x-y| > 8R } |\nabla_{x}G( x, y)|^2 \dx \biggr)^{\frac {1}{2}}> \lambda \biggr\} \lesssim  \lambda^{-\bar p} R^{-\frac{ d \bar p}{ 2} + d + (\bar p - p)}.
\end{align}
For any $0 < R \leq 1$ fixed, we may cover the set $\{|y| < 1\}$ by $N \lesssim R^{-d}$ balls of radius $R$; Hence, \eqref{est.gradient.jensen} and the second estimate in \eqref{est.C} yield that for all $\lambda > 0$
\begin{align}\label{inequality.covering.2}
&\pcap\biggl\{ |y| < 1 \, \colon \,  \M_0\biggl(\int_{ |x-y| > 8R } |G( x, y)|^{\frac{2d}{d-2}} \dx\biggr)^{\frac{d-2}{2d}} >  \lambda \biggr \} \lesssim \lambda^{- p} R^{-\frac{ d p}{ 2} },\\
&\pcap\biggl\{ |y| < 1 \, \colon \,  \M_0\biggl(\int_{ |x-y| > 8R } |\nabla_{x}G( x, y)|^2 \dx\biggr)^{\frac 1 2} >  \lambda \biggr \} \lesssim \lambda^{-\bar p} R^{-\frac{ d \bar p}{ 2} + (\bar p - p)}.\label{inequality.covering}
\end{align}
 Inequalities \eqref{est.00} immediately follow from these estimates if we choose $R=\frac 1 8$ and send $\lambda \uparrow +\infty$. 
 
 \medskip

We now derive \eqref{est.00b} for $\nabla_x G$ from \eqref{inequality.covering}. We begin by smuggling $R^\alpha$ into the left-hand side of \eqref{inequality.covering} and redefining $R^\alpha \lambda$ as $\lambda$ so that
\begin{align*}
\pcap\biggl\{ |y| < 1 \, \colon \, \M_0 \biggl( R^{2\alpha} \int_{ |x-y| > 8R}|\nabla_x G(x,y)|^2 \dx \biggr)^\oh > \lambda \biggr\}
\lesssim \lambda^{-\bar p} R^{-\frac{ d p}{ 2} + (\bar p - p)  +\alpha \bar p}.
\end{align*}
Since
\[ \int_{8R < |x-y|< 16R}|x-y|^{2\alpha} |\nabla_x G(x,y)|^2 \dx \lesssim R^{2\alpha}\int_{ |x-y|> 8R} |\nabla_x G(x,y)|^2 \dx,\]
we conclude that also
\begin{align}\label{est.A}
\pcap\biggl\{ |y| < 1 \, \colon \, \M_0 \biggl( \int_{8R< |x-y|< 16R}|x-y|^{2\alpha}|\nabla_x G(x,y)|^2 \dx \biggr)^\oh > \lambda \biggr\} \lesssim \lambda^{-\bar p} R^{-\frac{ d \bar p}{ 2} + (\bar p - p)+ \bar p \alpha}.
\end{align}

\smallskip

We now define
\begin{align}\label{set.1}
A&:= \biggl\{ |y| < 1 \, \colon \, \M_0 \biggl( \int_{|x-y|< 1}|x-y|^{2\alpha}|\nabla_x G(x,y)|^2 \dx \biggr)^\oh > \lambda \biggr\}
\end{align}
and, given a sequence of weights
\begin{align}
\{ &\omega_n \}_{n\in \mathbb{N}} \subset \R_+ \ \ \text{ such that $\sum_n \omega_n \leq 1$},\label{weights}
\end{align}
the sets
\begin{align}
A_n:= \biggl\{ |y| < 1 \, \colon \, \M_0 \biggl( & \int_{ 2^{-n} < |x-y|< 2^{-n+1}}\hspace{-0.3cm}|x-y|^{2\alpha}|\nabla_x G(x,y)|^2 \dx \biggr)^\oh > \omega_n \lambda \biggr\}\label{set.2}
\end{align}
for $n \in\N$. We claim that 
\begin{align}\label{est.B}
A \subset \bigcup_n A_n.
\end{align}
This can be easily seen by proving the  complementary statement  $\bigcap_n A_n^c \subset A^c$: Indeed, if
\[  \M_0 \biggl(  \int_{ 2^{-n} < |x-y|< 2^{-n+1}}|x-y|^{2\alpha}|\nabla_x G(x,y)|^2 \dx \biggr)^\oh \leq \omega_n \lambda \ \ \ \ \text{ for all $n\in \N$,}\]  
 then the inclusion of the sequence spaces $\ell^1 \subset \ell^2$, the subadditivity of the operator $\M_0$, and assumption \eqref{weights} yield
\begin{align*}
 \M_0 \biggl(& \int_{|x-y|< 1}|x-y|^{2\alpha}|\nabla_x G(x,y)|^2 \dx \biggr)^\oh\\
 & \leq \sum_{n\in \N}\M_0 \biggl(  \int_{ 2^{-n} < |x-y|< 2^{-n+1}}|x-y|^{2\alpha}|\nabla_x G(x,y)|^2 \dx \biggr)^\oh \leq  \lambda \sum_{n\in \N} \omega_n \stackrel{\eqref{weights}}{\leq} \lambda.
\end{align*}
We thus established  \eqref{est.B}.

\medskip

By \eqref{est.B}, we get
\[ \pcap ( A ) \leq \sum_n \pcap ( A_n ), \]
and, recalling definition \eqref{set.2}, we use estimate \eqref{est.A} with $\lambda$ and $R$ substituted by $\omega_n\lambda$ and $2^{-(n+3)}$ to bound
\begin{align*}
\pcap ( A ) \lesssim \lambda^{-\bar p}\sum_n \omega_n^{-\bar p} 2^{-n(-\frac{ d \bar p}{ 2} + (\bar p - p)+ \bar p \alpha)}.
\end{align*}
Choosing in \eqref{weights} $\omega_n= \frac{6}{(\pi n)^2}$, the sum on the right-hand side converges provided that $\alpha > \frac d 2 - \frac{\bar p -p}{\bar p}$ . 
Since we assumed $p < \bar p$, there exists $\alpha < \frac d 2$ such that
\begin{align*}
\pcap ( A ) \lesssim  \lambda^{-\bar p}.
\end{align*}
By definition \eqref{set.1}, sending $\lambda \uparrow +\infty$, we recover \eqref{est.00b} for $\nabla_x G$.

\medskip

The argument for $G$ in \eqref{est.00b} follows along the same lines as the one for $\nabla_xG$, as we shall argue now. In fact, since by H\"older's inequality with exponents $\frac{d}{d-2}$ and $\frac{d}{2}$ we may bound
$$
\biggl( \int_{ 8R < |x-y| < 16R } |G( x, y)|^2 \dx  \biggr)^{\frac 1 2} \lesssim R \biggl( \int_{ |x-y| > 8R } |G( x, y)|^{\frac{2d}{d-2}} \dx \biggr)^{\frac{d-2}{2d}},
$$
inequality \eqref{inequality.covering.2} yields
\begin{align*}
\pcap\biggl\{ |y-z| < 1 \, \colon \, R^{-1} \M_0\biggl( \int_{ 8R < |x-y| < 16R } |G( x, y)|^2 \dx \biggr)^\oh> \lambda \biggr\} \lesssim  \lambda^{- p} R^{-\frac{ d p}{ 2} },
\end{align*}
and thus also
\begin{align*}
\pcap\biggl\{ |y-z| < 1 \, \colon \, \M_0\biggl( \int_{ 8R < |x-y| < 16R } |G( x, y)|^2 \dx \biggr)^\oh> \lambda \biggr\} \lesssim  \lambda^{- p} R^{-\frac{ d p}{ 2} + p}.
\end{align*}
From this inequality we conclude  \eqref{est.00b} for $G$ as we did in the case of $\nabla_x G$ from \eqref{inequality.covering}. This concludes the proof of Step 3.

 \bigskip

 {\bf Step 4. Maximal function estimate.} We now prove \eqref{est.C} and begin with the first estimate. Without loss of generality, we focus on the case $z=0$. 
 For any $R > 0$ and $y \in \Rd$, let
\begin{align}\label{def.F}
 F_R(y) := \biggl( \int_{ |x| > 4R } |\nabla_{x}G( x, y)|^2 \dx \biggr)^{\frac{\bar p}{2p}}.
\end{align}
We first claim that it suffices to show that for every $R > 0$
\begin{align}\label{W1pF}
 \int_{|y| < R} \bigl(R^{-p} |F_R(y)|^p + |\nabla F_R(y)|^p\bigr) \dy  \lesssim R^{-\frac{d\bar p}{2} + d + (\bar p - p)}.
\end{align}
Indeed, if $\eta_R$ is a smooth cut-off function for $\{|y| < \frac R 2\}$ in $\{|y| < R\}$,  then by \eqref{W1pF} the product $\eta_{R} F_R$ satisfies 
$$
\| \eta_{R}F_R \|_{W^{1,p}(\Rd)}^p \stackrel{\eqref{W1pF}}{\lesssim} R^{-\frac{d\bar p}{2} + d + (\bar p - p)}.
$$
We thus apply the maximal function estimate \cite[Inequality (3.1)]{kinnunen_maximal}  to $\eta_{R}F_R$ and infer that
\begin{align*}
\pcap\biggl\{ |y| < R \, \colon \, \M(\eta_{R} F_R)(y) > \lambda   \biggr\} \lesssim \lambda^{-p} R^{-\frac{d\bar p}{2} + d + (\bar p - p)},
\end{align*}
where $\mathcal{M}$ is defined in \eqref{max.function}. Since by the assumption on $\eta_{R}$ and the definitions \eqref{max.function}$\&$\eqref{limsup} we have 
$$
\M_0 F_R \leq \M(\eta_{R}F_R) \ \ \ \ \text{ on $\{|y| < \frac R 2\}$,}
$$
we infer that
\begin{align}\label{est.CC}
 \pcap\biggl\{ |y| < \frac  R 2 \, \colon \, \M_0 F_R(y) > \lambda   \biggr\} \lesssim \lambda^{-p} R^{-\frac{d\bar p}{2} + d + (\bar p - p)}.
\end{align}
Furthermore, since $|y| < \frac R 2$ implies $\{ |x-y| > 8R \} \subset \{ |x| > 4R \}$ so that
\begin{align*}
\bigl( \int_{ |x-y| > 8R } |\nabla_{x}G( x, y)|^2 \, \dx \bigr)^{\frac{\bar p}{ 2p}} \leq F_R(y),
\end{align*}
 we conclude \eqref{est.C} for $\nabla_x G$ from \eqref{est.CC}.
 
\bigskip
 
To complete the argument for the first line in \eqref{est.C} it remains to prove \eqref{W1pF}:  The main ingredient for this are inequalities \eqref{grad.decay}$\&$\eqref{mixed.decay} which, by redefining $R$ as $2R$ and setting $z=0$, we rewrite as
\begin{align*}
\int_{|x | > 4R} \int_{|y| < 2R } \bigl( R^{-2}| \nabla_x G(x, y)|^2 + | \nabla_y\nabla_x G(x, y)|^2 \bigr) \dy  \dx \lesssim 1.
\end{align*}
Since by \eqref{symmetry.G} and \eqref{mixed.equation} the vector field $\nabla_x G( x, \cdot)$ is $a$-harmonic in $\{ |y| < 2R \}$
 for almost every $x$ with $|x | >  4R$, we may apply \eqref{Gehring2} of Step 1 to upgrade the previous estimate to
\begin{align*}
\int_{|x| > 4R} \biggl( \int_{|y| < R } R^{-p}| \nabla_x G(x, y)|^{\bar p} +  | \nabla_y\nabla_x G(x, y)|^{\bar p} \dy \biggl)^{\frac{ 2 }{\bar p}} \dx &\lesssim R^{-d  + \frac{2 d}{ \bar p}}.
\end{align*}
Since $\bar p \geq 2$, by Minkowski's inequality this in turn yields
\begin{align}
 \int_{|y| < R }\biggl(\int_{|x| > 4R}  | \nabla_y\nabla_x G(x, y)|^{2} \dx\biggl)^{\frac{\bar p }{2}}  \dy  &\lesssim R^{-\frac{d\bar p}{2}  + d},\label{Meyers.Green.1}\\
 \int_{|y| < R } \biggl(\int_{|x| > 4R}| \nabla_x G(x, y)|^{2} \dx\biggr)^{\frac {\bar p}{2}} \dy  &\lesssim R^{-\frac{d\bar p}{2}  + d + \bar p}.\label{Meyers.Green.2}
\end{align}

\smallskip

Differentiating \eqref{def.F} in $y$, the chain rule and Cauchy-Schwarz's inequality yield
\begin{equation}\label{der.y.F}
 |\nabla F_R(y)| \lesssim \biggl( \int_{|x| > 4R} |\nabla_x G(x, y)|^2 \d x \biggr)^{\frac{\bar p-p}{2p}} \biggl(\int_{|x| > 4R} |\nabla_x\nabla_y G(x,y)|^2 \d x\biggr)^{\frac 1 2}
\end{equation}
 By H\"older's inequality with exponents $\frac{\bar p}{\bar p - p}$  and  $\frac{\bar p}{p} > 1$, 
 \begin{align*}
 \int_{ |y| < R} &|\nabla F_R(y)|^p \dy \\
 &\leq \biggl(  \int_{ |y| < R } \bigl( \int_{|x| > 4R}  |\nabla_x G(x, y)|^2 \d x \bigr)^{\frac{\bar p}{2}} \dy \biggr)^{1- \frac{p}{\bar p}}  \biggl( \int_{ |y| < R }  \bigl(\int_{|x| > 4R} |\nabla_y\nabla_x G(x,y)|^2 \dx \bigr)^{\frac{\bar p}{ 2}} \dy \biggr)^{\frac{p}{\bar p}}.\notag
 \end{align*}
Inserting  \eqref{Meyers.Green.1}$\&$\eqref{Meyers.Green.2}, this implies \eqref{W1pF} for $\nabla F_R$. The $F_R$-part of  \eqref{W1pF} is immediate from \eqref{Meyers.Green.2}, which in view of \eqref{def.F}, assumes the desired form of
\begin{align*}
 \int_{ |y| < R } | F_R(y)|^p \dy  \lesssim   R^{-\frac{d\overline{p}}{2}+d+\overline{p}}. 
\end{align*}
This concludes the proof of \eqref{W1pF} and thus of the first estimate of \eqref{est.C}.

\medskip

The second estimate in \eqref{est.C} follows by a similar argument. We define for $R>0$ and $y\in \Rd$
\begin{align}\label{def.F.r.2}
F_R(y) : = \biggl( \int_{ |x| > 8R} |G(x,y)|^{\frac{2d}{d-2}} \dx \biggr)^{\frac{d-2}{2d}}.
\end{align}
Similarly to \eqref{der.y.F}, we may differentiate $F_R$ in the $y$-variable, apply H\"older's inequality, this time with exponents $\frac{2d}{d+2}$ and $\frac{2d}{d-2}$, and bound 
\begin{equation}\label{def.F.222}
|\nabla F_R(y)| \leq \biggl( \int_{ |x| > 8R} |\nabla_yG(x,y)|^{\frac{2d}{d-2} }\dx \biggr)^{\frac{d-2}{2d}}.
\end{equation}
Since by \eqref{equation.Green} and \eqref{mixed.equation} we have that for almost every $y \in \Rd$
 $$
 G( \cdot, y) \in Y^{1,2}_x(\{ |x- y|> 1\}), \ \ \nabla_y G(\cdot, y) \in Y^{1,2}_x(\{|x- y|>1 ; \Rd),
 $$
definition \eqref{def.F.r.2}, estimate \eqref{def.F.222} and Sobolev's inequality in the exterior domain $\{ |x| > 8R \}$ imply that 
\begin{align*}
\int_{|y| < R}|F_R(y)|^p \dy &\lesssim \int_{|y| < R}  \biggl( \int_{ |x| > 8R} |\nabla_x G(x,y)|^{2}\dx \biggr)^{\frac{p}{2}} \dy,\\
\int_{|y|< R} |\nabla F_R(y)|^p \dy &\lesssim \int_{|y| < R}  \biggl( \int_{ |x| > 8R} |\nabla_x\nabla_yG(x,y)|^{2}\dx \biggr)^{\frac{p}{2}} \dy .
\end{align*}
Another application of H\"older's inequality in $\{|y| < R \}$ with exponents $\frac{\bar p}{p}$ and $\frac{\bar p}{\bar p -p}$ further yields
\begin{align*}
\int_{|y| < R}|F_R(y)|^p \dy \lesssim R^{d(1- \frac{p}{\bar p})} \biggl(  \int_{|y| < R} \bigl( \int_{|x|> 8R} |\nabla_x G(x,y)|^2 \dx \bigr)^{\frac{\bar p }{2}} \dy \biggr)^{\frac{p}{\bar p}} \stackrel{\eqref{Meyers.Green.2}}{\lesssim} R^{-\frac{dp}{2} + d +p},\\
\int_{|y| < R}|\nabla F_R(y)|^p \dy \lesssim  R^{d(1- \frac{p}{\bar p})} \biggl( \int_{|y| < R}\bigl( \int_{|x|> 8R} |\nabla_x \nabla_y G(x,y)|^2 \dx \bigr)^{\frac{\bar p }{2}} \dy \biggr)^{\frac{p}{\bar p}}\stackrel{\eqref{Meyers.Green.1}}{\lesssim} R^{-\frac{dp}{2} + d},
\end{align*}
so that
\begin{align*}
 \int_{|y| < R} \bigl(R^{-p} |F_R(y)|^p + |\nabla F_R(y)|^p\bigr) \dy  \lesssim R^{-\frac{d p}{2} + d}.
\end{align*}
This is the analogue of \eqref{W1pF} for this definition of $F_R$. We may now pass from this to the second inequality in \eqref{est.C} as is done above to prove the first inequality in \eqref{est.C} from  \eqref{W1pF}. This concludes the proof of Step 4.

\bigskip

 {\bf Step 5. Construction of $G^*(a;  \cdot, \cdot)$.} \, Wrapping up Steps 2-4, we have that $G$ and $\nabla_x G$ satisfy \eqref{est.C} and therefore also \eqref{est.00}$\&$\eqref{est.00b} and \eqref{est.0}$\&$\eqref{est.0a} with an exponent $p$ that may be chosen strictly bigger than $2$. Equipped with these estimates, we now proceed to prove the existence of $G^*( \cdot, y)$ for $y$ outside an exceptional set $\Sigma$ satisfying (a) in the statement of Theorem \ref{t.main}.

\smallskip

For a test function $\zeta \in C^\infty_0(\Rd )$, we consider the function
\begin{align*}
u(y)= \int \zeta(x) G(x, y)  \dx.
\end{align*}
By the representation formula \eqref{representation.formula},  $u\in Y^{1,2}(\Rd)$ and  solves
\begin{align}\label{def.u}
-\nabla \cdot a\nabla u = \zeta \ \ \ \ \ \text{ in $\Rd$.}
\end{align}
Since $\zeta \in C^\infty_0(\Rd )$, Gehring's estimate \eqref{Gehring2} implies in particular that $u \in W^{1,p}_{loc}(\Rd)$. 
By Lebesgue's theorem for Sobolev functions \cite[Theorem 2.55]{ZiemerMaly}, we infer that the limit
\begin{align*}
\lim_{\delta \downarrow 0}& \fint_{|\tilde y -y| < \delta}u(\tilde y) \, d\tilde y= \lim_{\delta \downarrow 0} \fint_{|\tilde y -y| < \delta} \int \zeta(x) G( x, \tilde y) \dx \, d\tilde y
\end{align*}
exists as an element of $\R^m$ for all $y \in \Rd$ outside a set of zero $p$-capacity. Select a countable subset $\{ \zeta_n \}_{n \in \N} \subset C^\infty_0( \Rd)$ dense with respect to the $C^1$-topology. Hence, there exists a set $\tilde\Sigma$
with $\pcap( \tilde \Sigma)=0$ such that
\begin{align}\label{Leb.points}
\text{ $ \forall \, y \in \Rd \setminus \tilde\Sigma,$ $\, \forall n\in \N$}\ \ \lim_{\delta \downarrow 0} \fint_{|\tilde y -y| < \delta} \int \zeta_n(x) G(x, \tilde y)\dx \, d\tilde y \ \ \ \text{exists.}
\end{align}

\smallskip

Let $\Sigma_1$ and $\Sigma_2$ be the $p$-capacity zero sets  of \eqref{est.0}, \eqref{est.0a} in Step 2 and define
$$
\Sigma:= \tilde\Sigma \cup \Sigma_1 \cup \Sigma_2.
$$
With this definition, $\Sigma$ satisfies (a) of Theorem \ref{t.main}. By \eqref{est.00} and \eqref{est.00b} of Step 2, we remark that for every $y \notin \Sigma$, it also holds
\begin{align}\label{est.grad.G.max}
\mathcal{M}_0(\int \min\{ |x- y|^{2\alpha}, 1\} |\nabla_x G(x, y)|^2 \dx ) < +\infty,
\end{align}
for the exponent $\alpha < \frac d 2$ of Step 2.
By Jensen's inequality, for every $R > 0$ and $y \in \Rd$ we have
\begin{align*}
\limsup_{\delta \downarrow 0}\bigl\| \fint_{|\tilde y-y| < \delta } G( \cdot , \tilde y) \, d\tilde y& \bigr\|_{W^{1,q}_x(\{|x -y| < R\} )}\leq \M_0 \| G( \cdot , y)\|_{W^{1,q}_x(\{|x-y| < R \})},\\
\limsup_{\delta \downarrow 0}\bigl\| \fint_{|\tilde y-y| < \delta } G( \cdot , \tilde y) \, d\tilde y& \bigr\|_{Y^{1,2}_x(\{|x -y| > R\} )}\leq \M_0 \| G( \cdot , y)\|_{Y^{1,2}_x(\{|x-y| > R \} )}
\end{align*}
so that by \eqref{est.0}$\&$\eqref{est.0a} and weak compactness, we infer that for every $y\in \Rd \setminus \Sigma$ there exists a subsequence $\delta_k \downarrow 0$ (a priori depending on $y$) and a limit $G^*(a; \cdot, y) $ such that for all $R >0$
\begin{align}\label{def.ext.G}
\fint_{|\tilde y -y| < \delta_k} G( \cdot , \tilde y)\, d\tilde y \rightharpoonup G^*( \cdot , y)
\end{align}
in  $W^{1,q}_{x}(\{|x - y| < R \}) \cap Y^{1,2}_x(\{  |x -y | > R \})$. Moreover, inequality \eqref{est.grad.G.max} and weak lower-semicontinuity also yield
\begin{align}\label{est.grad.G}
\int \min\{ |x-y|^{2\alpha} , 1\} |\nabla_x G^*( x, y)|^2 \dx < +\infty.
\end{align}
We now show that \eqref{def.ext.G} holds for the entire family $\delta \downarrow 0$: Let us assume that this were not the case, i.e. that there exist two sequences $\{ \delta_k^{(1)}\}_k, \{\delta_k^{(2)}\}_{k}$ along which we obtain two different limits $G^{(1)}(\cdot, y), G^{(2)}(\cdot, y)$  in  \eqref{def.ext.G}. Appealing to \eqref{Leb.points}, to Fubini's theorem to exchange the order of the integrals, and to \eqref{def.ext.G} we infer that for every $n\in \N$
\begin{align*}
 \int \zeta_n(x) G^{(1)}( x, y) \dx = \int \zeta_n(x) G^{(2)}( x, y) \dx.
\end{align*}
Since the subset $\{ \zeta_n \}_{n\in \N}$ is chosen to be dense, we conclude that $G^{(1)}(x, y) = G^{(2)}(x, y)$ for almost every $x \in \Rd$. 

\smallskip

For every point $y$ outside $\Sigma$, we thus constructed $G^*( \cdot , y)$ which, by \eqref{def.ext.G}$\&$\eqref{est.grad.G}, satisfies (b) and the last inequality in (c) of Theorem \ref{t.main}. To conclude the proof of Theorem \ref{t.main}, it thus remains to show that $G^*( \cdot , y)$ solves equation \eqref{equation.Green}. Since $G( \cdot , \tilde y)$ solves equation \eqref{equation.Green} for almost every $\tilde y \in \Rd$, for every $\zeta \in C^\infty_0(\Rd)$, $y \in \Rd$ and $\delta > 0$ we have
\begin{align*}
\fint_{|\tilde y - y| < \delta} \int \nabla \zeta(x) \cdot a(x)\nabla_x G(x, \tilde y) \dx \d \tilde y =\fint_{|\tilde y- y| < \delta} \zeta(\tilde y) \d \tilde y,
\end{align*}
so that, by Fubini's theorem, 
\begin{align}\label{two.var.sol}
 \int \nabla \zeta(x) \cdot a(x) \fint_{|\tilde y- y| < \delta}\nabla_x G(x,\tilde y)\d \tilde y  \dx  = \fint_{|\tilde y- y| < \delta} \zeta(\tilde y) \d \tilde y.
\end{align}
Taking the limit $\delta \downarrow 0$, the assumption on $\zeta$, the boundedness \eqref{a.assumptions} of $a$ and \eqref{def.ext.G} yield that for all
$y \in \Rd \backslash \Sigma$ it holds
\begin{align}\label{two.var.sol}
 \int \nabla \zeta(x) \cdot a(x) \nabla_x G^*(x, y) \dx = \zeta( y).
\end{align}
Since $\zeta \in C^\infty_0(\Rd)$ is arbitrary, we conclude that $G^*( \cdot , y)$ solves equation \eqref{equation.Green}.

\smallskip

\noindent $\qed$

\bigskip

\noindent {\it Proof of Corollary \ref{cor1}.\,} Let $\Sigma$ be as in the statement of Theorem \ref{t.main} and $\eta$ be a smooth cut-off function for $\{|x| < 4 \}$ in $\{|x | < 2 \}$. With no loss of generality we may assume that each $u$ satisfies $\int_{|x| < 4} u= 0$.  The function $\eta u$ solves
\begin{align}\label{cut.harmonic}
 -\nabla \cdot a \nabla (\eta u) = \nabla \cdot g + f \ \ \ \ \text{in $\Rd$}
\end{align}
with 
\[ g := -u a \nabla \eta, \ \ \ f:= -\nabla \eta \cdot a \nabla u.\]
Both $g$ and $f$ are supported in $\{ 2 < |x| < 4 \}$ and by the definition of $\eta$, the second inequality in \eqref{a.assumptions}, the bound on the Dirichlet energy of $u$ and Poincar\'e's inequality, they satisfy
\begin{align}\label{bddrhs}
\int |g(x)|^2 \dx + \int |f(x)|^2 \dx  \lesssim 1.  
\end{align}

Furthermore, the representation formula \eqref{representation.formula},  \eqref{symmetry.G} and Theorem \ref{t.main} imply that for every $y \notin \Sigma$ with $|y|< 2 $ we may define as representative
\begin{align*}
u( y) = \int_{2< |x| <4 } \hspace{-0.3cm} g(x) \cdot \nabla_x G^*( x , y)\dx + \int_{2 < |x| < 4 }\hspace{-0.3cm} f(x) G^*( x, y) \dx.
\end{align*}

\smallskip

To make our notation leaner, we define
\begin{align}\label{splitting.representation}
v(y):= \int g(x) \cdot \nabla_x G^*( x, y) \dx, \ \ \ w(y):= \int f(x) G^*( x,y) \dx
\end{align}
and prove the statement of the corollary for  $v$. The function $w$ may be treated analogously.

\smallskip

We adapt the proof of \cite[Theorem 4.19]{EG} to show that there exist a sequence of sets $\{ B_j \}_{j \in \N} \subset \{ |y| < 1\}$ having
\begin{align}\label{p.cap.small}
\pcap\bigl\{ B_j \bigr\} < \frac {1}{2^{j}}
\end{align}
and moduli of continuity $\omega_j : \R_+ \to \R_+$ such that the following holds: For every $v$ constructed as in \eqref{splitting.representation}, there exists a sequence $\{ v_j \}_{j \in \N}$ satisfying for all $j \in \N$
\begin{align}\label{small.on.bj}
\sup_{\{|y| < 1\} \backslash B_j} |v(y) - v_j(y)| < \frac{1}{2^j}
\end{align}
and
\begin{align}\label{unif.cont}
|v_j(y) - v_j(\tilde y)| \leq \omega_j( |y - \tilde y|) \ \ \ \ \forall y, \tilde y \ \text{s.t.} \ |y | < 1, |\tilde y|< 1. 
\end{align}

\smallskip

From this, the corollary follows easily: For each $\eps > 0$ fixed, let $j_0 \in \N$ such that
$$
\pcap\bigl\{\bigcup_{j \geq j_0}B_j \bigr\} \stackrel{\eqref{p.cap.small}}{< }\frac \eps 2.
$$
By definition of capacity (see, for instance, \cite[Theorem 4.15 (i)]{EG}) we may find an open set $U^\eps \supset \bigcup_{j \geq j_0}B_j$ having $\pcap(U^\eps) < \eps$. We prove that on $\{ |y| < 1\}\backslash U^\eps$, the vector fields $v$ in \eqref{splitting.representation} are equicontinuous and that the continuity is uniform in $\{ |y| < 1\}\backslash U^\eps$. More precisely, this means proving that for each $\kappa >0$, there exists $\delta=\delta(\kappa) > 0$ such that for all $v$ as in \eqref{splitting.representation} and all $y, \tilde y \notin U^\eps$ and such that $|y| < 1$, $|\tilde y| < 1$ and $|y - \tilde y| < \delta$ we may bound
\begin{align*}
|v(y)- v(\tilde y) |< \kappa.
\end{align*}
By the triangle inequality, \eqref{small.on.bj} and the definition of $U^\eps$, we know indeed that if we fix $j \geq j_0$ such that $2^{-j} < \frac{\kappa}{3}$, then
\begin{align*}
|v(y) - v(\tilde y)| \leq |v(y) - v_{j}(y)| +  |v(\tilde y) - v_{j}(\tilde y)| +  |v^j(\tilde y) - v_{j}(y)| &\leq \frac{2}{3}\kappa + |v_j(\tilde y) - v_{j}(y)|\\
&  \stackrel{\eqref{unif.cont}}{\leq} \frac{2}{3}\kappa + \omega_j(|y -\tilde y|).
\end{align*}
It thus remains to pick $\delta$ such that the last term on the right-hand-side is smaller than $\frac \kappa 3$. This concludes the statement of the corollary.

\bigskip

We now show \eqref{p.cap.small}, \eqref{small.on.bj} and \eqref{unif.cont}. To do so, we begin by observing that, if $p> 2$ is as in Theorem \ref{t.main}, then by \eqref{Meyers.Green.1} in the proof of Theorem \ref{t.main}, we infer that
 \begin{align*}
 \int_{|y|< \frac 3 2}\bigl( \int_{2< |x| < 4} |\nabla_x \nabla_y G^*( x,y)|^2 \dx \bigr)^{\frac{p}{2}} \dy \lesssim 1.
 \end{align*}
 Similarly, this time using  \eqref{Meyers.Green.2}, we have that
 \begin{align*}
 \int_{|y|< \frac 3 2}\bigl( \int_{2< |x| < 4}& |\nabla_x G^*( x,y)|^2 \dx \bigr)^{\frac{p}{2}} \dy \lesssim 1.
 \end{align*}
 By standard approximation arguments adapted to Banach-valued functions (see e.g. \cite[Corollary 1.4.37]{CazenaveHaraux_BanachValued}), we may find a sequence $\{F_j \}_{j\in \N}$ of continuous maps 
 $$F_j: \bigl\{ |y| < \frac 3 2 \bigr\} \to L^2 \bigl(\{ 2 < |x| < 4 \} ; \R^{d} \bigr)$$
  such that for each $j \in \N$ we have
\begin{equation}\label{EG1}
\begin{aligned}
\int_{|y| < \frac 3 2} \bigl( \int_{2<|x| < 4}& |\nabla_x G^*( x,y) - F_j(x,y)|^2 \dx \bigr)^{\frac p 2} \dy \\
&+ \int_{|y| < \frac 3 2} \bigl( \int_{2<|x| < 4} |\nabla_y \nabla_x G^*( x,y) - \nabla_y F_j(x,y)|^2 \dx \bigr)^{\frac p 2}  \dy \leq \frac{1}{2^{(p+1)j}}.
\end{aligned}
\end{equation}

\bigskip

Let $\Sigma$ be the exceptional set of Theorem \ref{t.main} and let $\mathcal{M}_0$ be the maximal function operator (see \eqref{limsup}). We claim that
\begin{equation}\label{def.Bj}
\begin{aligned}
B_j:= \biggl\{ |y| < 1 \, \colon \, \ \ \mathcal{M}_0\bigl( \int_{2 < |x| < 4} | \nabla_x G^*(a; x, y) -  F_j(x, y)|^2& \dx \bigr)^{\frac 1 2} > \frac{1}{2^j}\biggr\} \cup \Sigma
\end{aligned}
\end{equation}
and $\omega_j$ such that for all $R> 0$
\begin{equation}\label{def.omega.j}
\omega_j(R) := 4 \sup \biggl\{  \biggl(\int_{2< |x| < 4} |  F_j(x,y) - F_j(x, \tilde y)|^2 \dx\biggr)^{\frac 12} \, | \, |y| \leq 1, \, |\tilde y| \leq 1, \,  |y - \tilde y| < R \biggr\},
\end{equation}
satisfy \eqref{p.cap.small}, \eqref{small.on.bj} and \eqref{unif.cont}
provided that for every $v$ we choose as approximating sequence  
\begin{align}\label{def.v.ij}
v_j(y):= \int g(x) \cdot F_j(x,y) \dx.
\end{align}
Here, $g$ is the vector field in the definition \eqref{splitting.representation} of $v$. We stress that, since each $F_j$ is continuous in $\{|y| \leq 1\}$ with values in $L^2 \bigl(\{ 2 < |x| < 4 \} ; \R^{d} \bigr)$, the above function is a well-defined modulus of continuity by the Heine-Cantor theorem. Furthermore, appealing to Cauchy-Schwarz's inequality, \eqref{bddrhs}  and \eqref{def.omega.j}, definition \eqref{def.v.ij}  immediately imply that $\{ v_j \}_{j \in \N}$ satisfy \eqref{unif.cont}.

\smallskip

It remains to show \eqref{p.cap.small} and \eqref{small.on.bj}: Since by Theorem \eqref{t.main} we have $\pcap( \Sigma)=0$, we may argue for $B_j$ as done for \eqref{est.C} and obtain that
\begin{equation}
\begin{aligned}\label{cap.Bi}
\pcap(B_j) &\lesssim 2^{pj } \int_{|y| < \frac 3 2} \bigl( \int_{2 < |x| < 4} | \nabla_x G^*( x,y) -  F_j(x,y)|^2 \dx \bigr)^{\frac p 2} \dy\\
&\quad\quad + 2^{pj } \int_{|y| < \frac 3 2}  \bigl( \int_{2 < |x| < 4} |\nabla_y  \nabla_x G^*( x,y) -  \nabla_y F_j(x,y)|^2 \dx \bigr)^{\frac p 2} \dy \stackrel{\eqref{EG1}}{\lesssim} 2^{-j},
\end{aligned}
\end{equation}
i.e. inequality \eqref{p.cap.small}. It remains to show that $\{ v_j \}_{j\in \N}$ defined in \eqref{def.v.ij} satisfies \eqref{small.on.bj}: For every $y \notin B_j \cup \Sigma$ with $|y|< 1$, we use the definition \eqref{splitting.representation} of $v$ and (b) of Theorem \ref{t.main} to rewrite
\begin{align*}
|v_j(y) - v(y) | =  \limsup_{r\downarrow 0} | v_j(y) - \fint_{|\tilde y - y| < r} v(\tilde y)\, \d\tilde y|.
\end{align*}
By the triangle inequality, we bound
\begin{align*}
|v_j(y) - v(y) | \leq   \limsup_{r\downarrow 0} | v_j(y) - \fint_{|\tilde y - y| < r} v_j(\tilde y)\, \d\tilde y| +   \limsup_{r\downarrow 0} |\fint_{|\tilde y - y| < r}( v(\tilde y) -v_j(\tilde y)) \, \d\tilde y| .
\end{align*}
By \eqref{unif.cont}, the first limit supremum on the right-hand side is zero. Hence, we have that 
\begin{equation}
\begin{aligned}\label{limsupsum}
| v_j(y) - v(y) |\leq  \limsup_{r\downarrow 0} \fint_{|\tilde y - y| < r} | v_j(\tilde y) - v(\tilde y)|\, \d\tilde y.
\end{aligned}
\end{equation}
Furthermore, the definitions of $v$ and $v_j$, Cauchy-Schwarz's inequality and the definition of $B_j$ together with \eqref{bddrhs} allow us to bound
\begin{align*}
&\limsup_{r \downarrow 0}\fint_{|\tilde y - y| < r } | v_j(\tilde y) -v(\tilde y)| \d\tilde y\\
&  \quad \lesssim \limsup_{r \downarrow 0} \fint_{|\tilde y -y|< r}\bigl( \int_{2< |x| < 4} | \nabla_x G^*( x,\tilde y) - F_j(x,\tilde y)|^2 \dx \bigr)^{\frac{1}{2}} \d\tilde y \lesssim  2^{-j}.
\end{align*} 
Inserting this into \eqref{limsupsum} we conclude \eqref{small.on.bj}.  The proof of Corollary \ref{cor1} is complete.

\smallskip

\noindent $\qed$

\section*{Acnowledgements}
The authors thank Jonas Lenz, University of Mainz, for pointing out an error in an earlier version of the paper.


\begin{thebibliography}{10}

\bibitem{CazenaveHaraux_BanachValued}
T.~Cazenave, A.~Haraux, and Y. Martel, \emph{An introduction to semilinear
  evolution equations}, Oxford Lecture Series in Mathematics and Its
  Applications 13, Clarendon Press; Oxford University Press, 1998.

\bibitem{ConlonGiuntiOtto}
J.~Conlon, A.~Giunti, and F.~Otto, \emph{Green function for elliptic systems:
  Delmotte-{D}euschel bounds}, Calc. of Var. and PDEs \textbf{56} (2017),
  no.~6.

\bibitem{DeGiorgiCounterexample}
E.~De~Giorgi, \emph{Un esempio di estremali discontinue per un problema
  variazionale di tipo ellittico}, Boll. Un. Mat. Ital.  \textbf{1} (1968), no. 4,
  135--137.

\bibitem{EG}
L.C. Evans and R.F. Gariepy, \emph{Measure theory and fine properties of
  functions}, CRC Press, 1992.

\bibitem{giaquinta}
M.~Giaquinta, \emph{Multiple integrals in the calculus of variations and
  nonlinear elliptic systems}, Annals of Mathematics Studies, vol. 105,
  Princeton University Press, Princeton, NJ, 1983.

\bibitem{JohnMalyStara}
O.~John, J. Mal\`y, and J.~Star\`a, \emph{Nowhere continuous solutions to
  elliptic systems}, Comment.Math.Univ.Carolinae \textbf{30} (1989), no.~1,
  33--43.

\bibitem{kinnunen_maximal}
J.~Kinnunen, \emph{The Hardy-Littlewood maximal function of a Sobolev
  function}, Israel Journal of Mathematics \textbf{100} (1997), 117--124.

\bibitem{ZiemerMaly}
J.~Mal\`y and W.P. Ziemer, \emph{Fine regularity of solutions of elliptic
  partial differential equations}, Providence, R.I. : American Mathematical
  Society, 1997.

\bibitem{MeyersLp}
N.G. Meyers, \emph{An $L^p$-estimate for the gradient of solutions of second
  order elliptic divergence equations}, Ann.Sc. Norm. Sup. Pisa \textbf{17}
  (1963), no.~3, 189--206.

\bibitem{Soucek}
J.~Soucek, \emph{Singular solution to linear elliptic systems},
  Comment.Math.Univ.Carolinae \textbf{25} (1987), 273--281.

\end{thebibliography}
\end{document}